\documentclass[final,3p]{elsarticle}
\usepackage{lineno,hyperref}
\usepackage{amsmath}
\usepackage{amssymb}
\usepackage[linesnumbered,ruled,vlined]{algorithm2e}
\usepackage{caption}
\usepackage{mathtools}
\usepackage{changes}
\usepackage{multirow}
\usepackage{verbatim}
\usepackage{mathrsfs}
\usepackage{graphicx}
\usepackage{subcaption}
\usepackage{amsmath,amsthm,bm,mathrsfs}
\theoremstyle{plain}

\theoremstyle{definition}

\theoremstyle{remark}

\modulolinenumbers[5]

\journal{ArXiv.org}

\begin{document}

\begin{frontmatter}

\title{A Low-rank Interpolatory Projection Algorithm for Solving Large-scale T-Sylvester Equations}

\author[uz]{Umair~Zulfiqar\corref{mycorrespondingauthor}}
\cortext[mycorrespondingauthor]{Corresponding author}
\ead{umair@yangtzeu.edu.cn}
\address[uz]{School of Electronic Information and Electrical Engineering, Yangtze University, Jingzhou, Hubei, 434023, China}
\begin{abstract}
This paper considers large-scale T-Sylvester equations of the form $AX - X^\top E^\top + B_1B_2^\top = 0$, which admit a low-rank solution. It is shown that when the unique solution of the T-Sylvester equation is low-rank, the problem naturally reduces to a tangential interpolation problem via oblique projection. The specific interpolation points and tangential directions needed to obtain the low-rank solution are not known a priori, thus requiring an iterative approach. An iterative interpolatory projection algorithm is proposed based on these interpolation conditions, which iteratively refines the interpolation data as the projection matrices expand in the number of columns. Numerical examples demonstrate that the proposed algorithm converges with projection matrices having significantly fewer columns compared to existing Krylov-subspace-based projection methods, confirming the superiority of the proposed algorithm over existing approaches.
\end{abstract}


\end{frontmatter}

\section{Introduction}
This paper considers the following T-Sylvester equation:
\begin{align}
AX-X^\top E^\top +B_1B_2^\top =0,\label{t_sylv}
\end{align}
where $A\in\mathbb{R}^{n\times n}$, $E\in\mathbb{R}^{n\times n}$, $B_1\in\mathbb{R}^{n\times m}$, $B_2\in\mathbb{R}^{n\times m}$, and $X\in\mathbb{R}^{n\times n}$.

Theoretical work on this equation dates to 1962 (Taussky and Wielandt \cite{taussky1962matrix}, for \(E = -A\)), with general solvability conditions established by Wimmer in the 1990s \cite{wimmer1994roth}. Recently, equation \eqref{t_sylv} has gained traction \cite{chiang2012sylvester,de2011consistency,de2011solution,dopico2016projection,garcia2013matrix} due to its link to palindromic eigenvalue problems (\(G + \lambda G^\top\)). Specifically, it is required for the first-order perturbation expansion of associated deflating subspaces \cite{byers2006structured}, i.e., Newton methods must repeatedly solve large-scale T-Sylvester equations \cite{chatelin1984simultaneous,demmel1987three}. It also serves as an auxiliary subproblem in structure-preserving QR algorithms for these eigenvalue problems \cite{kressner2009implicit}. For \(E = \pm A\), applications include Hamiltonian systems \cite{braden1998equations}, time-varying singular value decompositions (SVDs) \cite{baumann2003singular}, and quadratic inverse eigenvalue problems \cite{yuan2009direct}. For small-to-moderate matrices, modified Bartels–Stewart algorithms \cite{bartels1972algorithm} are typically used \cite{chiang2012sylvester,de2011consistency,vorontsov2011numerical}. Alternatively, \eqref{t_sylv} can be cast as an \(n^2 \times n^2\) linear system for standard iterative solvers \cite{wang2015computational}, though memory constraints restrict this approach to modest dimensions.

It is assumed that \(m \ll n\); under this assumption, the matrix \(X\) is numerically low-rank \cite{dopico2016projection}. Furthermore, both \(A\) and \(E\) are assumed to be invertible. Let \(\mu_1, \dots, \mu_n\) be the eigenvalues of \(E^{-1}A\). These eigenvalues satisfy the following conditions:
\begin{enumerate}
  \item \(\{\mu_1, \dots, \mu_n\} \setminus \{1\}\) is reciprocal‑free, i.e., for any pair \(\mu_i, \mu_j\) with neither equal to \(1\), we have \(\mu_i \mu_j \neq 1\).
  \item If \(1 \in \{\mu_1, \dots, \mu_n\}\), its algebraic multiplicity is \(1\); in other words, the eigenvalue \(\mu = 1\) can appear at most once in the spectrum of \(E^{-1}A\).
\end{enumerate}
These assumptions ensure that the T‑Sylvester equation \eqref{t_sylv} has a unique solution. For the general case where \(A\) and \(E\) are not invertible, the conditions for a unique solution can be found in Theorem 2.2 of \cite{dopico2016projection}.

Assume that the matrices \(\hat{V}\in\mathbb{R}^{n\times r}\) and \(\hat{W}\in\mathbb{R}^{n\times r}\) satisfy \(\hat{W}^\top\hat{W}=I\) and \(\hat{V}^\top\hat{V}=I\), where \(r\ll n\). Next, compute the following reduced-order matrices:
\begin{align}
\hat{A}=\hat{W}^\top A\hat{V},\quad \hat{E}=\hat{W}^\top E\hat{V},\quad \hat{B}_1=\hat{W}^\top B_1,\quad \hat{B}_2=\hat{W}^\top B_2,
\end{align}
with \(\hat{A}\in\mathbb{R}^{r\times r}\), \(\hat{E}\in\mathbb{R}^{r\times r}\), \(\hat{B}_1\in\mathbb{R}^{r\times m}\), and \(\hat{B}_2\in\mathbb{R}^{r\times m}\). Throughout the paper, the matrices \(\hat{A}\) and \(\hat{E}\) are assumed to be invertible.

Let \(\hat{X}\) be the unique solution of the following projected T‑Sylvester equation:
\begin{align}
\hat{A}\hat{X}-\hat{X}^\top \hat{E}^\top +\hat{B}_1\hat{B}_2^\top =0.\label{proj_t_sylv}
\end{align}
A projection‑based approximation of \(X\) can then be obtained as \(X\approx \hat{V}\hat{X}\hat{W}^\top\) \cite{dopico2016projection}.

The residual \(R\) for this approximation is defined by
\begin{align}
A\hat{V}\hat{X}\hat{W}^\top-\hat{W}\hat{X}^\top\hat{V}^\top E^\top +B_1B_2^\top =R,
\end{align}
and it satisfies the condition \(\hat{W}^\top R\hat{W}=0\).

The accuracy of the approximation \(X \approx \hat{V}\hat{X}\hat{W}^\top\) depends on the choice of the projection matrices \(\hat{V}\) and \(\hat{W}\). Two such choices are proposed in \cite{dopico2016projection}, as described below.

The first choice is
\begin{align}
\text{span}\left\{ \big(E^{-1}A\big)^{i-1}E^{-1}\begin{bmatrix}\hat{B}_1&\hat{B}_2\end{bmatrix} \mid i = 1,\dots,k \right\} \subseteq \text{Ran}(\hat{V}),\quad
\text{span}\left\{E\hat{V}\right\} \subseteq \text{Ran}(\hat{W}).\label{bk}
\end{align}

The second choice is
\begin{align}
&\text{span}\left\{ \big(E^{-1}A\big)^{i-1}E^{-1}\begin{bmatrix}\hat{B}_1&\hat{B}_2\end{bmatrix}, \big(A^{-1}E\big)^{i-1}A^{-1}\begin{bmatrix}\hat{B}_1&\hat{B}_2\end{bmatrix}\mid i = 1,\dots,k \right\} \subseteq \text{Ran}(\hat{V}),\nonumber\\
&\text{span}\left\{E\hat{V}\right\} \subseteq \text{Ran}(\hat{W}).\label{ebk}
\end{align}

In the next section, we propose an alternative choice of projection matrices, motivated by rational interpolation.
\section{Main Work}
In this section, we consider the case where \(X\) has rank \(r\). We then examine the truncated SVD of \(X = \hat{V}\hat{\Sigma}\hat{W}^\top\). Specifically, we identify the \(r\)-dimensional subspace \(\mathcal{V}\), whose basis consists of the columns of \(\hat{V}\), and the \(r\)-dimensional subspace \(\mathcal{W}\), whose basis consists of the columns of \(\hat{W}\). It is shown that \((zE - A)^{-1}\begin{bmatrix} B_1 & B_2 \end{bmatrix}\) is projected onto \(\mathcal{V}\) along \(\mathcal{W}\), with \(\hat{\Sigma}\) representing the solution to the projected T-Sylvester equation \eqref{proj_t_sylv}. The tangential interpolation conditions that \(\hat{V}\) and \(\hat{W}\) enforce are identified. Finally, a rational interpolation‑based projection algorithm for computing a low‑rank approximation of \(X\) is proposed.
\subsection{Interpolation Conditions Associated with the Low-rank Solution of \eqref{t_sylv}}
Consider the SVD of $X$ as $X=V\Sigma W^T$. The matrices $V$, $\Sigma$, and $W$ are then partitioned as follows:
\begin{align}
X=V\Sigma W^T=\begin{bmatrix}\hat{V}&\hat{V}_{n-r}\end{bmatrix}\begin{bmatrix}\hat{\Sigma}&0\\0&\Sigma_{n-r}\end{bmatrix}\begin{bmatrix}\hat{W}^\top\\W_{n-r}^\top\end{bmatrix}.
\end{align}
Next, we compute the quantities
\begin{align}
W^\top AV=\begin{bmatrix}\hat{A}&A_{12}\\A_{21}&A_{22}\end{bmatrix},\quad  W^\top EV=\begin{bmatrix}\hat{E}&E_{12}\\E_{21}&E_{22}\end{bmatrix},\quad W^\top B_1=\begin{bmatrix}\hat{B}_1\\B_{12}\end{bmatrix}, \quad  W^\top B_2=\begin{bmatrix}\hat{B}_2\\B_{22}\end{bmatrix}.
\end{align}
Recall that $X$ satisfies the T‑Sylvester equation \eqref{t_sylv}:
\begin{align}
AV\Sigma W^T-W\Sigma V^\top E^\top +B_1B_2^\top=0.\label{eq9}
\end{align}
Post‑multiplying \eqref{eq9} by $W$ yields
\begin{align}
AV\Sigma-W\Sigma V^\top E^\top W +B_1B_2^\top W&=0\nonumber\\
A\begin{bmatrix}\hat{V}&\hat{V}_{n-r}\end{bmatrix}\begin{bmatrix}\hat{\Sigma}&0\\0&\Sigma_{n-r}\end{bmatrix}
-\begin{bmatrix}\hat{W}&\hat{W}_{n-r}\end{bmatrix}\begin{bmatrix}\hat{\Sigma}&0\\0&\Sigma_{n-r}\end{bmatrix} \begin{bmatrix}\hat{E}^\top&E_{21}^\top\\E_{12}^\top&E_{22}^\top\end{bmatrix} +B_1\begin{bmatrix}\hat{B}_2\\B_{22}\end{bmatrix}^\top&=0\nonumber\\
A\hat{V}\hat{\Sigma}-\hat{W}\hat{\Sigma}\hat{E}^\top+B_1\hat{B}_2^\top-\hat{W}_{n-r}\Sigma_{n-r}E_{12}^\top&=0.
\end{align}
Pre‑multiplying \eqref{eq9} by $W^\top$ gives
\begin{align}
W^\top AV\Sigma W^T-\Sigma V^\top E^\top +W^\top B_1B_2^\top&=0,\nonumber\\
W\Sigma V^\top A^\top W- EV\Sigma +B_2 B_1^\top W&=0,\nonumber\\
\begin{bmatrix}\hat{W}&\hat{W}_{n-r}\end{bmatrix}\begin{bmatrix}\hat{\Sigma}&0\\0&\Sigma_{n-r}\end{bmatrix} \begin{bmatrix}\hat{A}^\top&A_{21}^\top\\A_{12}^\top&A_{22}^\top\end{bmatrix}
- E\begin{bmatrix}\hat{V}&\hat{V}_{n-r}\end{bmatrix}\begin{bmatrix}\hat{\Sigma}&0\\0&\Sigma_{n-r}\end{bmatrix} +B_2 \begin{bmatrix}\hat{B}_1\\B_{12}\end{bmatrix}&=0,\nonumber\\
\hat{W}\hat{\Sigma}\hat{A}^\top-E\hat{V}\hat{\Sigma}+B_2\hat{B}_1^\top+\hat{W}_{n-r}\Sigma_{n-r}A_{12}^\top&=0.
\end{align}
Now suppose $X$ has low rank $r$; consequently $\Sigma_{n-r}=0$ and 
\[
X=V\Sigma W^\top=\hat{V}\hat{\Sigma}\hat{W}^\top+\hat{V}_{n-r}\Sigma_{n-r}\hat{W}_{n-r}^\top=\hat{V}\hat{\Sigma}\hat{W}^\top.
\]
One easily sees that $\hat{X}=\hat{\Sigma}$ satisfies the projected T‑Sylvester equation \eqref{proj_t_sylv}.

Moreover, because $\Sigma_{n-r}=0$, we obtain
\begin{align}
A\hat{V}\hat{\Sigma}-\hat{W}\hat{\Sigma}\hat{E}^\top+B_1\hat{B}_2^\top&=0,\label{eq12}\\
\hat{W}\hat{\Sigma}\hat{A}^\top-E\hat{V}\hat{\Sigma}+B_2\hat{B}_1^\top&=0.\label{eq13}
\end{align}
Post‑multiplying \eqref{eq12} and \eqref{eq13} by $\hat{E}^{-\top}$ and $\hat{A}^{-\top}$, respectively, leads to
\begin{align}
A\hat{V}\hat{\Sigma}\hat{E}^{-\top}-\hat{W}\hat{\Sigma}+B_1\hat{B}_2^\top\hat{E}^{-\top}&=0,\label{eq14}\\
\hat{W}\hat{\Sigma}-E\hat{V}\hat{\Sigma}\hat{A}^{-\top}+B_2\hat{B}_1^\top\hat{A}^{-\top}&=0.\label{eq15}
\end{align}
Adding \eqref{eq14} and \eqref{eq15} yields
\begin{align}
A\hat{V}\hat{\Sigma}\hat{E}^{-\top}-E\hat{V}\hat{\Sigma}\hat{A}^{-\top}+B_1\hat{B}_2^\top\hat{E}^{-\top}+B_2\hat{B}_1^\top\hat{A}^{-\top}&=0.\label{eq16}
\end{align}
Post‑multiplying \eqref{eq16} by $\hat{E}^\top\hat{\Sigma}^{-1}$ gives
\begin{align}
A\hat{V}-E\hat{V}S_v+\begin{bmatrix}B_1&B_2\end{bmatrix}L_v=0,
\end{align}
where
$S_v=\hat{\Sigma}\hat{A}^{-\top}\hat{E}^\top\hat{\Sigma}^{-1}$, and $L_v=\begin{bmatrix}\hat{B}_2^\top\hat{\Sigma}^{-1}\\\hat{B}_1^\top\hat{A}^{-\top}\hat{E}^\top\hat{\Sigma}^{-1}\end{bmatrix}$.

Pre‑multiplying \eqref{eq12} and \eqref{eq13} by $-A^{-1}$ and $-E^{-1}$, respectively, we obtain
\begin{align}
-\hat{V}\hat{\Sigma}+A^{-1}\hat{W}\hat{\Sigma}\hat{E}^\top-A^{-1}B_1\hat{B}_2^\top&=0,\label{eq18}\\
-E^{-1}\hat{W}\hat{\Sigma}\hat{A}^\top+\hat{V}\hat{\Sigma}-E^{-1}B_2\hat{B}_1^\top&=0.\label{eq19}
\end{align}
Adding \eqref{eq18} and \eqref{eq19} results in
\begin{align}
A^{-1}\hat{W}\hat{\Sigma}\hat{E}^\top-E^{-1}\hat{W}\hat{\Sigma}\hat{A}^\top-A^{-1}B_1\hat{B}_2^\top-E^{-1}B_2\hat{B}_1^\top&=0.\label{eq20}
\end{align}
Post‑multiplying \eqref{eq20} by $\hat{E}^{-\top}\hat{\Sigma}^{-1}$ produces
\begin{align}
A^{-1}\hat{W}-E^{-1}\hat{W}S_w+\begin{bmatrix}A^{-1}B_1&E^{-1}B_2\end{bmatrix}L_w&=0,
\end{align}
with $S_w=\hat{\Sigma}\hat{A}^\top\hat{E}^{-\top}\hat{\Sigma}^{-1}$ and $L_w=\begin{bmatrix}-\hat{B}_2^\top\hat{E}^{-\top}\hat{\Sigma}^{-1}\\-\hat{B}_1^\top\hat{E}^{-\top}\hat{\Sigma}^{-1}\end{bmatrix}$.

Now perform eigenvalue decompositions of $S_v$ and $S_w$:
\[S_v=T_{v}^{-1}\begin{bsmallmatrix}\sigma_1&\cdots&0\\\vdots&\ddots&\vdots\\0&\cdots&\sigma_r\end{bsmallmatrix}T_{v}\quad\text{and}\quad S_w=T_{w}^{-1}\begin{bsmallmatrix}\mu_1&\cdots&0\\\vdots&\ddots&\vdots\\0&\cdots&\mu_r\end{bsmallmatrix}T_{w},\]
and define $b_i$, $c_i$ by
\[\begin{bmatrix}b_{1,1}&\cdots&b_{1,r}\\b_{2,1}&\cdots&b_{2,r}\end{bmatrix}=L_vT_{v}^{-1}\quad \text{and}\quad \begin{bmatrix}c_{1,1}&\cdots&c_{1,r}\\c_{2,1}&\cdots&c_{2,r}\end{bmatrix}=L_wT_{w}^{-1}.\]
Then, relying on the link between Sylvester equations and rational interpolation established in \cite{gallivan2004sylvester}, the following tangential interpolation conditions hold:
\begin{align}
(\sigma_iE-A)^{-1}\begin{bmatrix}B_1&B_2\end{bmatrix}\begin{bmatrix}b_{1,i}\\b_{2,i}\end{bmatrix}&=\hat{V}\big(\sigma_i\hat{E}-\hat{A}\big)^{-1}\begin{bmatrix}\hat{B}_1&\hat{B}_2\end{bmatrix}\begin{bmatrix}b_{1,i}\\b_{2,i}\end{bmatrix},\\
(\mu_iE^{-1}-A^{-1})^{-1}\begin{bmatrix}A^{-1}B_1&E^{-1}B_2\end{bmatrix}\begin{bmatrix}c_{1,i}\\c_{2,i}\end{bmatrix}&=\hat{W}\Big(\mu_i\hat{W}^\top E^{-1}\hat{W}-\hat{W}^\top A^{-1}\hat{W}\Big)^{-1}\nonumber\\
&\hspace*{2cm}\begin{bmatrix}\hat{W}^\top(A^{-1}B_1)&\hat{W}^\top (E^{-1}B_2)\end{bmatrix}\begin{bmatrix}c_{1,i}\\c_{2,i}\end{bmatrix}.
\end{align}
Furthermore, the projection matrices $\hat{V}$ and $\hat{W}$ satisfy
\begin{align}
\text{span}\left\{ \big(\sigma_iE-A\big)^{-1}\begin{bmatrix}B_1&B_2\end{bmatrix}\begin{bmatrix}b_{1,i}\\b_{2,i}\end{bmatrix}\mid i = 1,\dots,r \right\} &\subseteq \text{Ran}(\hat{V}),\\
\text{span}\left\{ \big(\mu_iE^{-1}-A^{-1}\big)^{-1}\begin{bmatrix}A^{-1}B_1&E^{-1}B_2\end{bmatrix}\begin{bmatrix}c_{1,i}\\c_{2,i}\end{bmatrix}\mid i = 1,\dots,r \right\} &\subseteq \text{Ran}(\hat{W});
\end{align}cf. \cite{gallivan2004sylvester}.

Now consider the Sylvester equations
\begin{align}
A\hat{V}_1-E\hat{V}_1\hat{A}^{-\top}\hat{E}^\top+B_1\hat{B}_2^\top&=0,\label{eq26}\\
A\hat{V}_2-E\hat{V}_2\hat{A}^{-\top}\hat{E}^\top+B_2\hat{B}_1^\top\hat{A}^{-\top}\hat{E}^\top&=0,\label{eq27}
\end{align}
where $\hat{V}=\big(\hat{V}_1+\hat{V}_2\big)\hat{\Sigma}^{-1}$.

Next, consider the Sylvester equations
\begin{align}
A^{-1}\hat{W}_1-E^{-1}\hat{W}_1\hat{A}^\top\hat{E}^{-\top}-A^{-1}B_1\hat{B}_2^\top\hat{E}^{-\top}&=0,\label{eq28}\\
A^{-1}\hat{W}_2-E^{-1}\hat{W}_2\hat{A}^\top\hat{E}^{-\top}-E^{-1}B_2\hat{B}_1^\top\hat{E}^{-\top}&=0,\label{eq29}
\end{align}
where $\hat{W}=\big(\hat{W}_1+\hat{W}_2\big)\hat{\Sigma}^{-1}$.

Pre‑multiplying \eqref{eq28} by $A$ gives
\begin{align}
EE^{-1}\hat{W}_1-AE^{-1}\hat{W}_1\hat{A}^\top\hat{E}^{-\top}-B_1\hat{B}_2^\top\hat{E}^{-\top}&=0,\nonumber\\
AE^{-1}\hat{W}_1\hat{A}^\top\hat{E}^{-\top}-EE^{-1}\hat{W}_1+B_1\hat{B}_2^\top\hat{E}^{-\top}&=0,\nonumber\\
AE^{-1}\hat{W}_1\hat{A}^\top-EE^{-1}\hat{W}_1\hat{A}^\top\hat{A}^{-\top}\hat{E}^\top+B_1\hat{B}_2^\top&=0.\label{eq30}
\end{align}
Assuming unique solutions for the Sylvester equations \eqref{eq26} and \eqref{eq30}, one readily observes that $W_1=EV_1\hat{A}^{-\top}$.

Similarly, pre‑multiplying \eqref{eq29} by $E$ yields
\begin{align}
EA^{-1}\hat{W}_2-AA^{-1}\hat{W}_2\hat{A}^\top\hat{E}^{-\top}-B_2\hat{B}_1^\top\hat{E}^{-\top}&=0,\nonumber\\
AA^{-1}\hat{W}_2\hat{A}^\top\hat{E}^{-\top}-EA^{-1}\hat{W}_2+B_2\hat{B}_1^\top\hat{E}^{-\top}&=0,\nonumber\\
AA^{-1}\hat{W}_2\hat{E}^\top\hat{E}^{-\top}\hat{A}^\top-EA^{-1}\hat{W}_2\hat{E}^\top+B_2\hat{B}_1^\top&=0,\nonumber\\
AA^{-1}\hat{W}_2\hat{E}^\top-EA^{-1}\hat{W}_2\hat{E}^\top\hat{A}^{-\top}\hat{E}^{\top}+B_2\hat{B}_1^\top\hat{A}^{-\top}\hat{E}^{\top}&=0.\label{eq31}
\end{align}
If the Sylvester equations \eqref{eq27} and \eqref{eq31} have unique solutions, we similarly obtain $W_2=AV_2\hat{E}^{-\top}$.

Thus, once $\hat{V}_1$ and $\hat{V}_2$ are computed, $\hat{W}_1$ and $\hat{W}_2$ follow immediately. The matrices $\hat{V}_1$, $\hat{V}_2$, $\hat{W}_1$, $\hat{W}_2$ possess the following properties:
\begin{align}
\text{span}\left\{ \big(\sigma_iE-A\big)^{-1}B_1b_{1,i}\mid i = 1,\dots,r \right\} &\subseteq \text{Ran}(\hat{V}_1),\\
\text{span}\left\{ \big(\sigma_iE-A\big)^{-1}B_2b_{2,i}\mid i = 1,\dots,r \right\} &\subseteq \text{Ran}(\hat{V}_2),\\
\text{span}\left\{ E\hat{V}_1 \right\} &\subseteq \text{Ran}(\hat{W}_1),\\
\text{span}\left\{ A\hat{V}_2 \right\} &\subseteq \text{Ran}(\hat{W}_2).
\end{align}
We now rewrite the Sylvester equations \eqref{eq26} and \eqref{eq27} as
\begin{align}
A(\hat{V}_1\hat{E}^{-\top})-E(\hat{V}\hat{E}^{-\top})(\hat{E}^{-1}\hat{A})^{-\top}+B_1(\hat{E}^{-1}\hat{B}_2)^\top&=0,\\
A(\hat{V}_2\hat{E}^{-\top}\hat{A}^\top\hat{E}^{-\top})-E(\hat{V}\hat{E}^{-\top}\hat{A}^\top\hat{E}^{-\top})(\hat{E}^{-1}\hat{A})^{-\top}+B_2(\hat{E}^{-1}\hat{B}_1)^\top&=0.
\end{align}
Next, decompose $\hat{E}^{-1}\hat{A}$ via its eigenvalue decomposition:
\[\hat{E}^{-1}\hat{A}=T_l^{*}\begin{bsmallmatrix}\hat{\lambda}_1&\cdots&0\\\vdots&\ddots&\vdots\\0&\cdots&\lambda_r\end{bsmallmatrix}T_l^{-*},
\]
and define $l_{1,i}$, $l_{2,i}$ by
\[\begin{bmatrix}l_{1,1}\\\vdots\\l_{1,r}\end{bmatrix}=T_{l}^{-*}\big(\hat{E}^{-1}\hat{B}_1\big) \quad \text{and}\quad \begin{bmatrix}l_{2,1}\\\vdots\\l_{2,r}\end{bmatrix}=T_{l}^{-*}\big(\hat{E}^{-1}\hat{B}_2\big) .\]
One then observes that the projection matrices $\hat{V}_1$ and $\hat{V}_2$ satisfy
\begin{align}
\text{span}\left\{ \Big(\frac{1}{\hat{\lambda}_i}E-A\Big)^{-1}B_1l_{2,i}^\top\mid i = 1,\dots,r \right\} &\subseteq \text{Ran}(\hat{V}_1),\\
\text{span}\left\{ \Big(\frac{1}{\hat{\lambda}_i}E-A\Big)^{-1}B_2l_{1,i}^\top\mid i = 1,\dots,r \right\} &\subseteq \text{Ran}(\hat{V}_2).
\end{align}
These conditions resemble the rational interpolation conditions arising in model order reduction for discrete‑time systems. In the special case $m=1$ or $B_1=B_2$, $\hat{V}_1$ and $\hat{V}_2$ essentially satisfy a subset of the optimality conditions for the $\mathcal{H}_2$‑optimal model order reduction problem for discrete‑time systems; see \cite{bunse2010h2} for details.

In summary, when $X$ has low rank, the task of finding a low‑rank solution naturally transforms into a tangential‑interpolation‑based projection problem. The specific tangential interpolation conditions satisfied by the projected matrices, as well as the $r$‑dimensional subspaces $\mathcal{V}$ onto which the matrices of the T‑Sylvester equation \eqref{t_sylv} are projected, have been identified. The next subsection develops a rational‑interpolation‑based algorithm founded on these interpolation conditions.
\subsection{Rational Interpolation-based Low-rank Solver}
The interpolation points $\frac{1}{\hat{\lambda}_i}$ and the tangential directions $l_{1,i}$, $l_{2,i}$ needed to compute the projection matrices $\hat{V}_1$ and $\hat{V}_2$ depend on the reduced matrices $\hat{A}$, $\hat{E}$, $\hat{B}_1$, and $\hat{B}_2$ associated with the low‑rank solution of $X$, which are not known \textit{a priori}. This situation closely resembles $\mathcal{H}_2$‑optimal model order reduction for discrete‑time state‑space models, where the optimal interpolation points are the reciprocals of the poles of the optimal reduced‑order model and the tangential directions are the residues associated with those poles. Since the optimal reduced‑order model is unknown in advance, an iterative algorithm was proposed in \cite{bunse2010h2}, which starts with an initial guess of the interpolation points and tangential directions and repeatedly updates the interpolation data until convergence. In our setting, unlike in $\mathcal{H}_2$‑optimal model order reduction, $r$ is not fixed; instead, it grows until the norm of the residual $R$ falls below a desired tolerance. Nevertheless, the choice of interpolation data can be made in the same manner as in rational‑interpolation‑based $\mathcal{H}_2$‑optimal model order reduction for discrete‑time state‑space models.

Starting from an arbitrary initial guess of the interpolation point $\sigma_1$ and tangential directions $b_{1,1}$, $b_{2,1}$, we compute $\hat{V}_1$ and $\hat{V}_2$ as follows:
\begin{align}
&\text{span}\left\{(\sigma_1E-A)^{-1}B_1b_{1,1}\right\} \subseteq \text{Ran}(v_1),\quad \text{span}\left\{ (\sigma_1E-A)^{-1}B_2b_{2,1}\right\} \subseteq \text{Ran}(v_2).
\end{align}
We then set $\hat{V}=\mathrm{orth}\big(\begin{bmatrix}v_1&v_2\end{bmatrix}\big)$ and $\hat{W}=\mathrm{orth}\big(\begin{bmatrix}Ev_1&Av_2\end{bmatrix}\big)$, and compute the reduced‑order matrices $\hat{E}$, $\hat{A}$, $\hat{B}_1$, $\hat{B}_2$. Thereafter, the next interpolation points $\sigma_i$ and tangential directions $b_{1,i}$, $b_{2,i}$ are obtained as
$\sigma_i=\frac{1}{\hat{\lambda}_i}$, $b_{1,i}=l_{2,i}^\top$, and $b_{2,i}=l_{1,i}^\top$. The projection matrices $\hat{V}$ and $\hat{W}$ can be orthogonally expanded (e.g., using modified Gram–Schmidt) as
\[
\hat{V}=\mathrm{orth}\big(\begin{bmatrix}\hat{V}&v_1&v_2\end{bmatrix}\big)\quad \text{and} \quad \hat{W}=\mathrm{orth}\big(\begin{bmatrix}\hat{W}&Ev_1&Av_2\end{bmatrix}\big).
\]

As the columns of $\hat{V}$ and $\hat{W}$ increase, the number of eigenvalues of $\hat{E}^{-1}\hat{A}$ grows. It is well understood in the rational interpolation and eigenvalue literature that using the most controllable poles of the triplet $(\hat{E},\hat{A},[\hat{B}_1,\hat{B}_2])$ as interpolation points generally yields good accuracy in rational interpolation \cite{gugercin2008h_2} and also accelerates convergence in pole estimation for eigensolvers \cite{rommes2007methods,rommes2006,rommes2006efficient,mengi2022large}. Note that $\hat{E}^{-1}\hat{A}$ contains Ritz values of $E^{-1}A$, and as $\hat{V}$ and $\hat{W}$ expand in columns, it is expected to capture the poles of $E^{-1}A$. Although this is not our primary objective, the interpolation data can be generated automatically without user intervention, similarly to eigensolvers such as subspace‑accelerated dominant pole estimation (SADPA) \cite{rommes2006,rommes2006efficient}. The most controllable pole of $\hat{E}^{-1}\hat{A}$ is the one associated with the largest value of  
\[
\hat{\phi}_i = \frac{\|[l_{1,i},l_{2,i}]\|_2^2}{|\mathrm{Re}(\hat{\lambda}_i)|},
\]
as discussed in \cite{rommes2007methods,mengi2022large}. This criterion is widely used in eigensolvers to identify the most controllable poles of $\hat{E}^{-1}\hat{A}$.

If the number of columns of $\hat{V}$ and $\hat{W}$ becomes excessively large without a sufficient reduction in the residual $R$, making the eigenvalue decomposition of $\hat{E}^{-1}\hat{A}$ expensive, one can apply an implicit restart as in SADPA \cite{rommes2006,rommes2006efficient}. This essentially means retaining only the last few columns of $\hat{V}$ and $\hat{W}$ and discarding the earlier history. Moreover, to avoid repeatedly generating the same interpolation point in successive iterations, deflation can be applied, as is done in most sparse eigensolvers like SADPA \cite{rommes2006,rommes2006efficient}.

One advantage of tangential interpolation is that, regardless of the value of $m$, two columns are added to $\hat{V}$ and $\hat{W}$ each time they are orthogonally expanded. However, the accuracy of tangential interpolation depends on both appropriate interpolation points and tangential directions. If $m$ is small, one can instead use block interpolation and define $\hat{V}_1$, $\hat{V}_2$ as
\begin{align}
&\text{span}\left\{(\sigma_iE-A)^{-1}B_1\right\} \subseteq \text{Ran}(v_1^{(i)}),\quad \text{span}\left\{ (\sigma_iE-A)^{-1}B_2\right\} \subseteq \text{Ran}(v_2^{(i)}),
\end{align}
in which case $2m$ columns are added to $\hat{V}$ and $\hat{W}$ during each orthogonal expansion.

The pseudo‑code of the proposed low‑rank solver for large‑scale T‑Sylvester equations, named “A Low‑rank Interpolatory Projection Algorithm for T‑Sylvester Equations (LRIPA‑TSYLV)”, is presented in Algorithm \ref{alg}. The stopping criterion
\[
\frac{\|A\hat{V}\hat{X}\hat{W}^\top-\hat{W}\hat{X}^\top\hat{V}^\top E^\top +B_1B_2^\top\|_F}{(\|A\|_F+\|E^\top\|_F)\|\hat{X}\|_F-\|B_1B_2^\top\|_F} \leq \tau
\]
used in LRIPA‑TSYLV is identical to the one used in the projection algorithms presented in \cite{dopico2016projection}.
\begin{algorithm}[h]
\DontPrintSemicolon
\SetAlgoLined
\SetKwInOut{KwIn}{Inputs}
\SetKwInOut{KwOut}{Output}
\KwIn{Matrices of T-Sylvester Equation \eqref{t_sylv}: $(E,A,B_1,B_2)$; Type: ``\textit{block}'' or ``\textit{tangential}''; Tolerance: $\tau$.}
\KwOut{Low-rank approximation of $X$: $X\approx\hat{V}\hat{X}\hat{W}^\top$.}

$i = 1$, $\sigma_i = 0$, $b_{1,i} = \mathbf{1}_{m \times 1}$, $b_{2,i} = \mathbf{1}_{m \times 1}$, $B_{\perp} = \begin{bmatrix}B_1 & B_2\end{bmatrix}$, $\hat{V} = [\;]$, $\hat{W} = [\;]$.\;
\While{$\frac{\|A\hat{V}\hat{X}\hat{W}^\top-\hat{W}\hat{X}^\top\hat{V}^\top E^\top +B_1B_2^\top\|_F}{(\|A\|_F+\|E^\top\|_F)\|\hat{X}\|_F-\|B_1B_2^\top\|_F} > \tau$}{
    \eIf{\textnormal{Type} = ``\textit{block}''}{
        $\hat{B} = \begin{bmatrix}B_1 & B_2\end{bmatrix}$.\;
    }{
        $\hat{B} = \begin{bmatrix}B_1 b_{1,i} & B_2 b_{2,i}\end{bmatrix}$.\;
    }
    Solve for $\begin{bmatrix}v_1^{(i)} & v_2^{(i)}\end{bmatrix}$: $(\sigma_i E - A)\begin{bmatrix}v_1^{(i)} & v_2^{(i)}\end{bmatrix} = \hat{B}$.\;
    \eIf{$\mathrm{Im}(\sigma_i) = 0$}{
        $\hat{V} \gets \mathrm{orth}\Big(\begin{bmatrix}\hat{V} & v_1^{(i)} & v_2^{(i)}\end{bmatrix}\Big)$ and $\hat{W} \gets \mathrm{orth}\Big(\begin{bmatrix}\hat{W} & E v_1^{(i)} & A v_2^{(i)}\end{bmatrix}\Big)$.\;
    }{
        $\hat{V} \gets \mathrm{orth}\Big(\begin{bmatrix}\hat{V} & \mathrm{Re}(v_1^{(i)}) & \mathrm{Im}(v_1^{(i)}) & \mathrm{Re}(v_2^{(i)}) & \mathrm{Im}(v_2^{(i)})\end{bmatrix}\Big)$ and
        $\hat{W} \gets \mathrm{orth}\Big(\begin{bmatrix}\hat{W} & \mathrm{Re}(E v_1^{(i)}) & \mathrm{Im}(E v_1^{(i)}) & \mathrm{Re}(A v_2^{(i)}) & \mathrm{Im}(A v_2^{(i)})\end{bmatrix}\Big)$.\;
    }
    Project: $\hat{A} \gets \hat{W}^\top A \hat{V}$, $\hat{E} \gets \hat{W}^\top E \hat{V}$, $\hat{B}_1 \gets \hat{W}^\top B_1$, $\hat{B}_2 \gets \hat{W}^\top B_2$, $\hat{B}_{\perp} \gets \hat{W}^\top B_{\perp}$.\;
    Solve the projected T-Sylvester equation \eqref{proj_t_sylv} to compute $\hat{X}$.\;
    Deflate: $B_{\perp} \gets \begin{bmatrix}B_1 & B_2\end{bmatrix} - E \hat{V} \hat{E}^{-1} \begin{bmatrix}\hat{B}_1 & \hat{B}_2\end{bmatrix}$.\;
    Compute the eigenvalue decomposition of $\hat{E}^{-1}\hat{A}$: $\hat{E}^{-1}\hat{A} = \hat{T}\,\mathrm{diag}(\hat{\lambda}_1,\dots,\hat{\lambda}_r)\hat{T}^{-1}$.\;
    Set $r_{k} = \hat{T}^{-1}(k,:)\hat{E}^{-1}\hat{B}_{\perp}$ and sort the columns of $\hat{T}$ in descending order of $\hat{\phi}_k = \frac{\|r_{k}\|_2^2}{|\mathrm{Re}(\hat{\lambda}_k)|}$.\;
    Update: $i \gets i+1$, $\sigma_i \gets \frac{1}{\hat{T}^{-1}(1,:)\hat{E}^{-1}\hat{A}\hat{T}(:,1)}$, $b_{1,i} \gets \big(\hat{T}^{-1}(1,:)\hat{E}^{-1}\hat{B}_2\big)^\top$, $b_{2,i} \gets \big(\hat{T}^{-1}(1,:)\hat{E}^{-1}\hat{B}_1\big)^\top$.\;
}
\caption{LRIPA-TSYLV}
\label{alg}
\end{algorithm}
\section{Numerical Results}
In this section, we assess the numerical performance of LRIPA‑TSYLV and compare it with the Block Krylov method for T‑Sylvester equations (BK‑TSYLV, Algorithm 1 in \cite{dopico2016projection}) and the Extended block Krylov method for T‑Sylvester equations (EBK‑TSYLV, Algorithm 2 in \cite{dopico2016projection}), which respectively use the projection matrices given in \eqref{bk} and \eqref{ebk}. The test problems are taken from the numerical experiments section of \cite{dopico2016projection}. MATLAB codes for reproducing the results are publicly available at \cite{mycode}. The convergence tolerance $\tau$ is fixed at $10^{-10}$ for all examples. All computations are performed using MATLAB R2025b on a Windows 11 laptop equipped with 32 GB of RAM and an Intel(R) Core(TM) Ultra 9 285H processor running at 2.9 GHz.
\subsection{Example 1}
The matrices $A$ and $E$ are of size $10^4 \times 10^4$ and arise from finite‑difference discretizations on the domain $[0,1] \times [0,1]$ of the differential operators
\begin{align}
a(u) &= \bigl(-\exp(-xy)\, u_x\bigr)_x + \bigl(-\exp(xy)\, u_y\bigr)_y + 100\, x\, u_x + \gamma\, u\,,\nonumber\\
e(u) &= -u_{xx} - u_{yy} + 100\, x\, u_x\,,\nonumber
\end{align}
with $\gamma = 5 \cdot 10^4$. The vectors $B_1\in\mathbb{R}^{10^4\times 1}$ and $B_2\in\mathbb{R}^{10^4\times 1}$ are generated using MATLAB's \textit{rand} command and are then scaled by $10^4$ to match the magnitude of the entries of the matrices $A$ and $E$; see \cite{dopico2016projection}. The maximum number of columns for $\hat{V}$ and $\hat{W}$ is set to $200$.

\begin{figure}[!h]
  \centering
  \includegraphics[width=12cm]{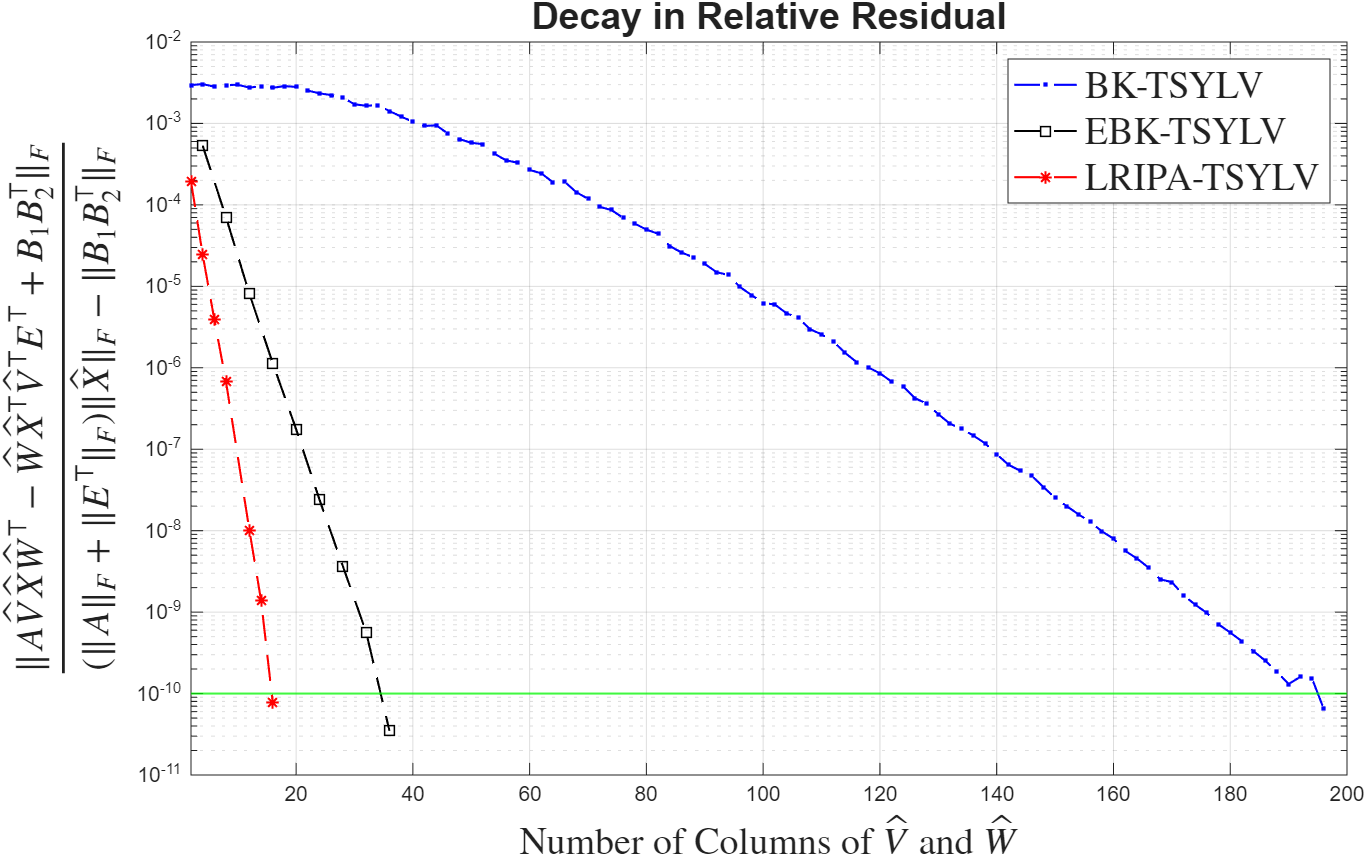}
  \caption{Decay in relative residual for \(X\)}\label{fig1}
\end{figure}

Since $m=1$ in this example, setting Type = ``\textit{block}'' or Type = ``\textit{tangential}'' in LRIPA‑TSYLV makes no difference. BK-TSYLV, EBK-TSYLV, and LRIPA‑TSYLV converged in $12.4452$ sec, $0.3012$ sec, and $0.2878$ sec, respectively. The relative residual's history is plotted in Figure \ref{fig1}. It can be seen that LRIPA‑TSYLV produced the most compact approximation of $X$.
\subsection{Example 2}
The matrices $A$ and $E$, each of size $10^4 \times 10^4$, are obtained from finite‑difference discretizations over the domain $[0,1] \times [0,1]$ of the differential operators
\begin{align}
a(u) &= (-\exp(-xy)\, u_x)_x + (-\exp(xy)\, u_y)_y + 100\, x\, u_x + \gamma\, u,\nonumber\\
e(u) &= -u_{xx} - u_{yy},\nonumber
\end{align}
where $\gamma = 5 \cdot 10^4$. The matrices $B_1\in\mathbb{R}^{10^4\times 2}$ and $B_2\in\mathbb{R}^{10^4\times 2}$ are generated using MATLAB's \textit{rand} command and then scaled by $10^4$ to match the magnitude of the entries of $A$ and $E$ \cite{dopico2016projection}. The maximum number of columns for $\hat{V}$ and $\hat{W}$ is set to $200$.

\begin{figure}[!h]
  \centering
  \includegraphics[width=12cm]{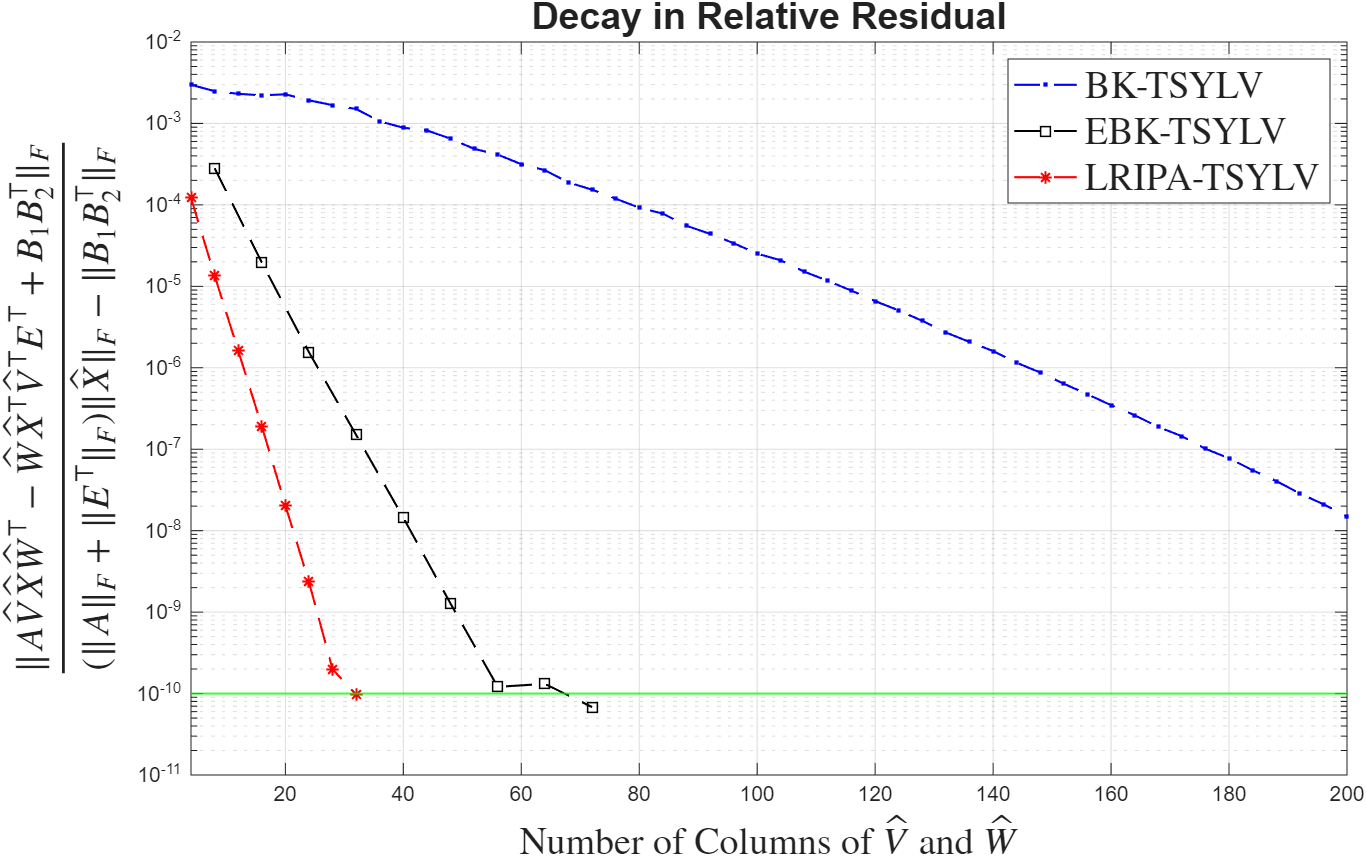}
  \caption{Decay in relative residual for \(X\)}\label{fig2}
\end{figure}

Since $m$ in this example is small, we set Type = ``\textit{block}'' in LRIPA‑TSYLV. EBK‑TSYLV and LRIPA‑TSYLV converged in $0.5321$ sec and $0.3934$ sec, respectively, while BK‑TSYLV did not converge and took $7.1516$ sec. The history of the relative residual is plotted in Figure \ref{fig2}. It can be seen that LRIPA‑TSYLV produced the most compact approximation of $X$.
\subsection{Example 3}
The matrices $A$ and $E$ are of size $4\cdot10^4 \times 4\cdot10^4$ and are obtained via finite‑difference discretizations on the domain $[0,1] \times [0,1]$ of the differential operators
\begin{align}
a(u) &= -u_{xx} - u_{yy} + y(1-x)u_x + \gamma u\,,\nonumber\\
e(u) &= -u_{xx} - u_{yy}\,,\nonumber
\end{align}
respectively, with $\gamma = 10^4$. The matrices $B_1\in\mathbb{R}^{4\cdot10^4\times 5}$ and $B_2\in\mathbb{R}^{4\cdot10^4\times 5}$ are generated using MATLAB's \textit{rand} command and then scaled by $10^4$ to match the magnitude of the entries of $A$ and $E$ \cite{dopico2016projection}. The maximum number of columns for $\hat{V}$ and $\hat{W}$ is set to $500$.

\begin{figure}[!h]
  \centering
  \includegraphics[width=12cm]{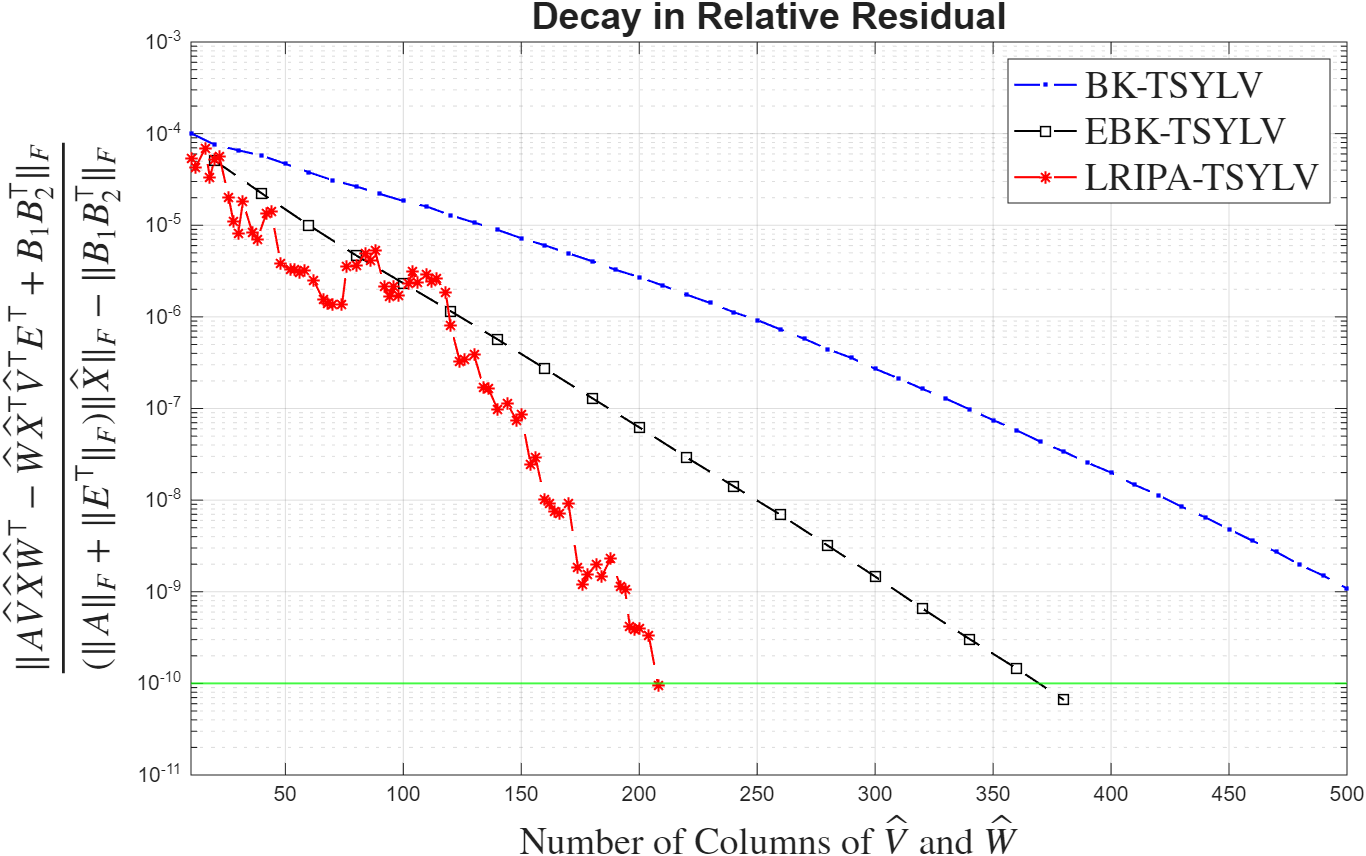}
  \caption{Decay in relative residual for \(X\)}\label{fig3}
\end{figure}

We set Type = ``\textit{tangential}'' in LRIPA‑TSYLV so that only two new columns are added to $\hat{V}$ and $\hat{W}$ in each iteration; the block version would add $10$ columns per iteration and cause $\hat{V}$ and $\hat{W}$ to grow rapidly. EBK‑TSYLV and LRIPA‑TSYLV converged, while BK‑TSYLV did not. The history of the relative residual is plotted in Figure \ref{fig3}. It can be seen that LRIPA‑TSYLV produced the most compact approximation of $X$. The computational time with and without residual computation is reported in Table \ref{tab}. EBK‑TSYLV took less time than LRIPA‑TSYLV despite producing a less compact approximation. This is because LRIPA‑TSYLV requires eigenvalue decompositions to generate interpolation data, whereas EBK‑TSYLV does not. Moreover, the LU factorizations for the linear solves in EBK‑TSYLV are reused, while in LRIPA‑TSYLV they cannot be reused since the interpolation point $\sigma_i$ changes at each iteration. These factors make LRIPA‑TSYLV computationally less efficient, even though it offers the advantage of producing a compact approximation $\hat{V}\hat{X}\hat{W}^\top$ that requires less memory to store.
\begin{table}[!h]
\centering
\caption{Elapsed Time Comparison}\label{tab}
\begin{tabular}{|l|l|}
\hline
Method                         & Elapsed Time (sec) \\ \hline
BK-TSYLV                       & 64.2642           \\ 
BK-TSYLV (without residual computation)    & 45.6595           \\ 
EBK-TSYLV                      & 18.2391            \\ 
BK-TSYLV (without residual computation)    & 12.9302            \\
LRIPA-TSYLV                    & 37.6109           \\
LRIPA-TSYLV (without residual computation) & 27.0858           \\ \hline
\end{tabular}
\end{table}
\section{Conclusion}
In this paper, we have addressed the efficient numerical solution of large‑scale T‑Sylvester equations that admit a low‑rank solution. We have shown that computing a low‑rank solution naturally reduces to a tangential interpolation problem. The required interpolation points and tangential directions depend on the reduced matrices, which are not known a priori. To overcome this difficulty, we have developed an iterative interpolatory projection algorithm, LRIPA‑TSYLV, that progressively builds the projection subspaces \(\hat{V}\) and \(\hat{W}\) by extracting interpolation data from the current reduced matrices and then expanding the subspaces accordingly. The proposed algorithm produces a compact approximation that typically requires far smaller dimensions than existing Krylov‑subspace‑based methods, leading to significantly lower storage requirements. By using tangential interpolation, the growth of the projection matrices is controlled—only two new columns are added per iteration. The algorithm automatically generates the interpolation points and tangential directions without user interference. LRIPA‑TSYLV is an efficient algorithm for computing low‑rank solutions of large‑scale T‑Sylvester equations.
\section*{Acknowledgment}
The author is grateful to Prof. Froilán M. Dopico for sharing his MATLAB implementation of Algorithm 3.1 from his paper \cite{de2011consistency}.

\end{document}